\documentclass[12pt, reqno]{amsart}
\usepackage{fullpage}
\usepackage{graphicx}
\usepackage{bbm}
\newtheorem{theorem}{Theorem}[subsection]
\theoremstyle{plain}

\newtheorem{corollary}{Corollary}[subsection]

\newtheorem{lemma}{Lemma}[subsection]

\newtheorem{problem}{Problem}
\newtheorem{proposition}{Proposition}[subsection]

\newcommand{\customqed}[1]{{\renewcommand{\qedsymbol}{#1}\qed}}
\newcommand{\varqed}{\customqed{\hbox{$\lrcorner$}}}

\theoremstyle{definition}

\newtheorem{definitionAux}{Definition}[subsection]
\newenvironment{definition}{\begin{definitionAux}}{\varqed\end{definitionAux}}

\theoremstyle{remark}

\newtheorem{remarkAux}{Remark}[subsection]
\newenvironment{remark}{\begin{remarkAux}}{\varqed\end{remarkAux}}

\newtheorem{exampleAux}{Example}[subsection]

\newtheorem{remarksAux}{Remarks}[subsection]
\newenvironment{remarks}{
\begin{remarksAux}\ \begin{enumerate}}{\end{enumerate}\varqed\end{remarksAux}}

\newtheorem{examplesAux}{Examples}[subsection]
\newenvironment{examples}{
\begin{examplesAux}\ \begin{enumerate}}{\end{enumerate}\varqed\end{examplesAux}}

\numberwithin{equation}{section}

\renewcommand{\S}{\mathcal{S}}
\newcommand{\N}{\mathbb{N}}
\newcommand{\uple}[1]{{\langle #1 \rangle}}

\newcommand{\R}{\mathbb{R}}
\newcommand{\Z}{\mathbb{Z}}
\newcommand{\norm}[1]{\left\Vert #1 \right\Vert}
\newcommand{\defin}[1]{\emph{\textbf{#1}}}

\renewcommand{\O}{\mathcal{O}}

\newcommand{\Q}{\mathbb{Q}}

\newcommand{\integral}[2]{\int\!{#1}\,\mathrm{d}{#2}}

\input{tcilatex}

\begin{document}
\date{March 12, 2009}
\title{Dynamics and abstract computability:
computing invariant measures.}
\author{Stefano Galatolo, Mathieu Hoyrup and Crist\'obal Rojas}

\begin{abstract}
We consider the question of computing invariant measures from an abstract
point of view. We work in a general framework (computable metric spaces,
computable measures and functions) where this problem can be posed precisely. We 
consider invariant measures as fixed points of the transfer operator and give general conditions under which the transfer operator is (sufficiently) computable. In this case, a general result ensures the computability of isolated fixed points and hence invariant measures (in given classes of ``regular'' measures).
This implies the computability of  many SRB measures. 

On the other hand, not all computable dynamical systems have a computable invariant measure. We exhibit two interesting examples of  computable dynamics, one having an SRB measure which is not computable and another having no computable invariant measure at all, showing some subtlety in this kind of problems.
\end{abstract}

\maketitle

\tableofcontents

\section{Introduction}

An important fact motivating the study of the statistical properties of
dynamical systems is that the pointwise long time prediction of a chaotic
system is not possible, while the estimation or forecasting of averages and
other long time statistical properties is sometimes possible. This often
corresponds in mathematical terms to computing invariant measures, or estimating
some of their properties.

Giving a precise meaning to the computation of a continuous object like a measure is not a completely obvious task and involves the definition of effective versions of several  concepts from mathematical analysis.

Our approach will be mainly based on the concept of computable metric space.
To give a first example, let us consider the set $\mathbb{R}$ of real
numbers. Beyond $\mathbb{Q}$ there are many other real numbers that can be
handled by algorithms: $\pi $ or $\sqrt{2}$ for instance can be approximated {\em at any given precision} (with rational numbers) by an
algorithm. Hence these numbers can be identified with the algorithm which is
able to calculate them (more precisely, with the string representing the program
which approximates it). This set of points is called the set of \emph{computable real numbers} and was introduced in the famous paper \cite{T36}.

This kind of construction can then be generalized to many other metric spaces,
considering a dense countable set that plays the same role as the rationals in
the above example. Then, \emph{computable} or \emph{recursive} counterparts
of many mathematical notions can be defined, and rigorous statements about
the algorithmic approximation of abstract objects can be made, also obtaining algorithmic versions of many classical theorems (see section 2). In particular, this general approach also gives the possibility to treat in a simple way measures spaces, define computable measures and computable functions between measure spaces (transfer operators), which will be the main theme of this paper.

The paper is devoted to the problem of computation of invariant measures in discrete time dynamical systems. By discrete time dynamical system we mean a system $(X,T)$ were $X$ is a metric space and $T:X\to X$ is a Borel measurable transformation. Here an invariant measure is a Borel measure $\mu$ on $X$ such that for each measurable set $A$ it holds $\mu(A)=\mu(T^{-1}(A))$. Such measures contain information on the statistical behavior of the system $(X,T)$ and on the possible behavior of averages of observables  along typical trajectories of the system. The map $T$ moreover induces a function $L_T:PM(X)\to PM(X)$, where $PM(X)$ is the set of Borel probability measures over $X$ and will be endowed with a suitable metric (for details see section \ref{meas}). 
$L_T$ is called the transfer operator associated to $T$ (basic results about this are reminded in section \ref{DS}).

Before entering into details about the computation of measures and invariant measures in particular, we remark that whatever we mean by \textquotedblleft approximating a measure by an algorithm\textquotedblright , there are only countably many \textquotedblleft measure approximating algorithms\textquotedblright\ whereas, in general, a dynamical system may have uncountably many invariant measures (usually an infinite dimensional set). So, most of them will not be algorithmically describable. 
This is not a serious problem because we can put our attention on the most \textquotedblleft meaningful\textquotedblright\ ones. An important part of the theory of dynamical systems is indeed devoted to the understanding of
\textquotedblleft physically\textquotedblright\ relevant invariant measures,
among these, SRB measures play an important role\footnote{Informally speaking, these are measures which represent the asymptotic statistical
behavior of \textquotedblleft many\textquotedblright\ (positive Lebesgue
measure) initial conditions, see section \ref{DS}}. These measures are good candidates to be
computed. The existence and uniqueness of SRB measures is a widely studied 
problem (see \cite{Y}), which has been solved for some important classes of
dynamical systems.

Let us precise the concept of computable measure. As mentioned before, the framework of computable analysis can be applied to abstract spaces as the space $PM(X)$. 
A measure $\mu $ is then \emph{computable} if it is a computable point of that measure space.  In this case there is an algorithm such that, for each rational $\varepsilon $ given as input, outputs a "finite" measure (a finite rational convex combination of Dirac measures supported on \textquotedblleft rational\textquotedblright\ points) which is $\varepsilon $-close to $\mu $.

In the literature, there are several works dealing with the problem of approximating invariant measures, more or less informally from the algorithmic point of view (see e.g. \cite{L}, \cite{H}, \cite{KMY}, \cite{PJ99}, \cite{Din93, Din94}). In these works the main technique consists in an adequate discretization of the problem.
 More precisely, in several of the above works the transfer operator associated to the dynamics (see subsection \ref{PF}) is approximated by a finite dimensional one and the problem is reduced to the  computation of the corresponding relevant eigenvectors (some convergence result then validates the quality of the approximation).

Another strategy to face the problem of computation of invariant measures consist in following the way the measure $\mu $ can be constructed and check that each step can be realized in an effective way. In some interesting examples we can obtain the SRB measure as limit of iterates of the Lesbegue measure
$\mu =\lim_{n\rightarrow \infty }L_{T}^{n}(m)$ where $m$ is the Lesbegue
measure and $L_{T}$ is the transfer operator associated to $T$. To prove computability of $\mu $ the main point is to recursively estimate the speed of convergence to the limit. This sometimes can be done using the decay of correlations (see \cite{GHR07} where computability of SRB measures in uniformly hyperbolic systems is proved in this way, see \cite{GP09} for general relations between convergence of measures and decay of correlations with a point of view similar to the one of the present paper).

Let us illustrate the main results of the paper. Informally speaking, a function $T:X\rightarrow X$ is said to be computable if its behavior can be described by some algorithm (for the precise definitions see sections \ref{CMS} and \ref{comp_meas}). In this case the pair $(X,T)$ is called a computable dynamical system. In this context, the general problem we are facing can be stated in the following terms:

\begin{problem}\label{probl}

\begin{enumerate}
\item[{\bf a)}]  Given a computable dynamical system $(X,T)$ does the set of invariant
measures contain computable points?

\item[{\bf b)}] Can they be found in an algorithmic way, starting from the description of the system?
\end{enumerate}
\end{problem}

We will see that, in general, even the above question a) does not always have a
positive answer. However, in many interesting situations, both of the above problems can be positively solved.

We will take a general point of view finding the interesting invariant
measure as a fixed point of the transfer operator, giving general
conditions ensuring its computability. The following theorem will be the
main tool (see Thm. \ref{comp_inv}).\newline

\noindent \textbf{Theorem A}\textit{\ Let $X$ be a computable metric space
and $T$ a function which is computable on $X\setminus D$. Let us consider
the dynamical system $(X,T)$. Suppose there is a recursively compact set of
probability measures $V\subset PM(X)$ such that for every $\mu \in V$, $\mu
(D)=0$ holds. Then every invariant measure isolated (in the weak topology) in $V$ is computable.}\newline

The precise meaning of computability on $X\setminus D$ will be given in
section \ref{CMS} however the intuitive meaning of the above proposition is
that: if the function $T$ is computable outside some singular set $D$ (the
discontinuity set for example) and we look for invariant measures in a set $V
$ of measures giving no weight to the set $D$ (some class of regular
measures e.g.) and in the set $V$ there is a unique invariant measure, then
this measure can be computed.

This will give as a consequence that the SRB measure is computable in many
examples of computable systems (uniquely ergodic systems, piecewise expanding maps in one
dimensions, systems having an indifferent fixed point and many other systems
having an unique absolutely continuous invariant measure, see Theorem \ref{comp_meas} and Prop. \ref{2} ).

Observe that any object which is \textquotedblleft
computable\textquotedblright\ in some way (as $T,V,\mu $ in the theorem)
admits a finite description (a finite program). Theorem A is actually \emph{%
uniform}: there is a \emph{single} algorithm which takes finite descriptions
of $T$ and $V$ and which, as soon as the hypothesis in Theorem A are satisfied
and $\mu $ is a unique invariant measure in $V$, outputs a finite
description of $\mu $ (see remark \ref{comprem} and the above item b) of
Problem \ref{probl}). Observe that the algorithm cannot decide whether the hypotheses are satisfied or not, but computes the measure whenever they are fulfilled.

After such general statements, one could conjecture that, in computable
dynamical systems, SRB measures are always computable. This is not true, and
reveals some subtlety about the general problem of computing an invariant
measure. In section \ref{lst} we will see that:\newline

\noindent \textbf{Examples}\textit{\ There exists a computable dynamical
system having no computable measure at all. Moreover, there exists a
computable dynamical system on the unit interval having an SRB measure which
is not computable.}\newline

The interest of the second example comes from the fact that any computable
map of the interval must have some computable invariant measure. The example
shows that important invariant measures can still be missed.

To further motivate these results, we finally remark that from a technical
point of view, computability of the considered measure is a requirement in
several results about relations between computation, probability, randomness
and pseudo-randomness (see e.g. \cite{LM08}, \cite{GHR08}, \cite{GHR07},\cite{GHR09c}).

\subsection{Plan of the paper}

In section 2 we give a  compact and self contained introduction to
the prerequisites about computable analysis which are necessary to work with
dynamical systems on metric spaces. In this section we also prove some
general statements about solutions of equations on metric spaces which will
be used to ``find'' the interesting invariant measures as fixed points of
the transfer operator (Theorem, \ref{computable_fixed}).

In section 3 we develop the computable treatment of the space of probability
measures on a given (computable) metric space. Some results of these initial sections
are new and should be of independent interest. Their usefulness is
demonstrated in the next sections.

In section 4 we start considering dynamical systems. A direct application of
the results of the previous sections allow us to establish general
assumptions under which the transfer operator is computable (on a suitable
subset, Theorem \ref{L_comp}).

We then use the framework and tools introduced before to face Problem \ref{probl}. We prove Theorem A
above (which also becomes a simple application of previous results) and show
how to apply it in order to prove the computability of many interesting
invariant measures.


In section 6 we construct the two counter-examples already announced.

\section{Preliminaries on algorithmic theory}
\subsection{Analysis and computation}
A way to approach several problems from mathematical analysis by
computational tools is to approximate the ``infinite'' mathematical objects
(elements of non countable sets, as real numbers or a functions ) involved in
the problem by some algorithm which constructs an approximating sequence of
``finite'' objects (rational numbers, polynomials with rational coefficients)
which are ``treatable'' by the computer. Usually, the algorithm has to
manipulate and decide questions about the various mathematical objects
involved, and convergence results should be provided in order to choose the
suitable level of accuracy for the finite approximation. The actual
implementation of the algorithm and the various decisions are, in most
cases, subjected to round off errors which can produce additional
approximation errors, wrong decisions or undecidable situations if the error
is not considered rigorously (how to decide $x\geq y$ when $x=y$?).
Sometimes, estimates (for these errors) can be obtained under suitable
conditions, but this is in general a further and often nontrivial task (see
e.g. \cite{Bla94}). In this paper we will work in a framework where the
algorithmic abilities of the computer to represent and manipulate infinite
mathematical objects are taken into account from the beginning. In this
framework (often referred to as Computable Analysis) one can rigorously
determine which objects can be algorithmically approximated at any given
accuracy (these will be called \emph{computable} objects), and which can not.

 
Here, the word \emph{computable} is used, but may be adapted to each
particular situation: for instance, ``computable'' functions from $\mathbb{N}$
to $\mathbb{N}$ are called \emph{recursive} functions, ``computable'' subsets
of $\mathbb{N}$ are called r.e sets, etc.

\subsection{Background from recursion theory}\label{background}

The starting point of recursion theory was to give a mathematical definition
making precise the intuitive notions of \emph{algorithmic} or \emph{%
effective procedure} on symbolic objects. Every mathematician has a more or
less clear intuition of what can be computed by algorithms: the
multiplication of natural numbers, the formal derivation of polynomials are
simple examples.

Several very different formalizations have been independently proposed (by
Post, Church, Kleene, Turing, Markov\dots) in the 30's, and have proved to
be equivalent: they compute the same functions from $\mathbb{N}$ to $\mathbb{%
N}$. This class of functions is now called the class of \emph{recursive
functions}. As an algorithm is allowed to run forever on an input, these
functions may be \emph{partial}, i.e.~not defined everywhere. The \emph{%
domain} of a recursive function is the set of inputs on which the algorithm
eventually halts. A recursive function whose domain is $\mathbb{N}$ is said
to be \emph{total}. For formal definitions see for example~\cite{Rog87}.

With this intuitive description it is more or less clear that there exists
an effective procedure to enumerate the class of all partial recursive
functions, associating to each of them its \emph{\textbf{G\"odel number}}.
Hence there exists a \emph{universal} recursive function $\varphi_u:\mathbb{N}
\rightarrow \mathbb{N}$ satisfying for all $e,n\in \mathbb{N}$, $\varphi_u({%
\langle e,n \rangle})=\varphi_e(n)$ where $e$ is the G\"odel number of $%
\varphi_e$ and ${\langle \cdot,\cdot \rangle}:\mathbb{N}^2 \rightarrow 
\mathbb{N}$ is some recursive bijection.

The notion of recursive function induces directly an important computability
notion on the class of subsets of $\mathbb{N}$: a set of natural numbers is
said to be \emph{\textbf{recursively enumerable}} (\emph{\textbf{r.e}} for
short) if it is the range of some partial recursive function. That is if
there exists an algorithm listing (or enumerating) the set. We denote by $%
E_e $ the r.e set associated to $\varphi_e$, namely: $E_e=\mathrm{range}%
(\varphi_e)=\{\varphi_u({\langle e,n \rangle}):n\in \mathbb{N}\}$, where $\varphi_{u}$ is the universal recursive function.

Let $(E_i)_{i\in\mathbb{N}}$ be a family of r.e subsets of $\mathbb{N}$. We
say that $E_i$ is r.e \emph{\textbf{uniformly in $\boldsymbol{i}$}} if there
is a single recursive function $\varphi$ such that $E_i=\{\varphi({\langle
i,n \rangle}):n\in\mathbb{N}\}$. Taking $\varphi=\varphi_u$ the universal
recursive function yields an enumeration $(E_i)_{i\in\mathbb{N}}$ of all the
r.e subsets of $\mathbb{N}$, such that $E_i$ is r.e uniformly in $i$.

More generally, once a computability notion has been defined for some class
of objects in the following form:

\begin{center}\label{computable_object}
An object $x$ is \emph{computable} if there is a (partial or total)  recursive function $\varphi $ which computes $x$ in some sense.
\end{center}

 A uniform version will
be implicitly defined and intensively used:

\begin{center}
Objects from a family $(x_i)_{i\in\mathbb{N}}$ of $X$ are \emph{uniformly}
computable if there is a single (total or partial) recursive function $%
\varphi$ such that $\varphi({\langle i,. \rangle}):\mathbb{N}\to\mathbb{N}$
computes $x_i$ for each $i$.
\end{center}


\subsection{From $\mathbb{N}$ to countable sets}

Strictly speaking, recursive functions only work on natural numbers, but
this can be extended to the objects (thought of as ``finite'' objects) of
any countable set, once a numbering of its elements has been chosen.

\begin{definition}
A \emph{\textbf{numbered set}} $\O $ is a countable set together with a
surjection $\nu_\O:\mathbb{N} \to \O $ called the \emph{\textbf{numbering}}.
We write $o_n$ for $\nu(n)$ and call $n$ a \emph{\textbf{name}} of $o_n$.
\end{definition}


The set $\mathbb{Q}$ of rational numbers can be injectively numbered $%
\mathbb{Q}=\{q_0,q_1,\ldots\}$ in an \emph{effective} way: the number $i$ of
a rational $a/b$ can be computed from $a$ and $b$, and vice versa. We fix
such a numbering. %
%


\begin{definition}
A subset $A$ of a numbered set $\O $ is \emph{\textbf{recursively enumerable
(r.e)}} if there is a r.e set $E\subseteq \mathbb{N}$ such that $%
A=\{o_n:n\in E\}$.
\end{definition}

Uniformity for r.e subsets of $\O $ is defined as uniformity for r.e subsets
of $\mathbb{N}$.


\subsection{Computability of reals}

The following notion was already introduced by Turing in \cite{T36}.

\begin{definition}
Let $x$ be a real number. We say that:

$\bullet$ $x$ is \emph{\textbf{lower semi-computable}} if the set $\{q\in%
\mathbbm{Q}:q<x\}$ is r.e.,

$\bullet$ $x$ is \emph{\textbf{upper semi-computable}} if the set $\{q\in%
\mathbbm{Q}:q>x\}$ is r.e.,

$\bullet$ $x$ is \emph{\textbf{computable}} if it is lower and upper
semi-computable.
\end{definition}

The following classical characterization may be more intuitive: a real
number is computable if and only if there exists a recursive function $%
\varphi$ computing a sequence of rational numbers converging exponentially
fast to $x$, that is $|q_{\varphi(i)}-x|<2^{-i}$, for all $i$. We remark that as there exists subsets of integers which are recursively enumerable but not recursive (see \cite{Rog87}), there also exists semi-computable numbers which are not computable.

In the following section we will see how these notions can be generalized to
separable metric spaces, which inherit the computable
structure of $\mathbbm{R}$ via the metric.


\subsection{Computable metric spaces\label{CMS}}

In this section we introduce the basic tools of computable analysis on
metric spaces. Most of the results of this section and several of the following one have been already  obtained by Weihrauch, Brattka, Presser and others in the framework of ``Type-2 theory of Effectivity'', which is based in the notion of ``representation'' (infinite binary codes) of mathematical objects. A standard reference book on this approach to  Computable Analysis is \cite{Wei00}, and a specific paper on computability of subsets of metric spaces is \cite{BraPre03}.
Our approach to Computable Analysis only uses the notion of recursive function (see subsection \ref{background}). It is intended to emphasize the fact that computability notions are just the ``effective'' versions of classical ones. In this way we obtain a theory syntactically familiar to most mathematicians and computability results can be proved in a transparent and compact way.
 
A computable metric space is a metric space with a dense
numbered set such that the distance on this set is algorithmically compatible
with the numbering (distances between numbered points can be uniformly
computed up to arbitrary precision). From this point of view the real line (with euclidean distance) has a natural structure of computable metric space, whit the rationals as a numbered set.

\begin{definition}
A \emph{\textbf{computable metric space}} (CMS) is a triple $\mathcal{X}%
=(X,d,\S )$, where

$\bullet$ $(X,d)$ is a separable complete metric space,

$\bullet$ $\S =(s_i)_{i \in \mathbbm{N}}$ is a dense subset of $X$ (the
numbered set of \emph{\textbf{ideal points}}),

$\bullet$ The real numbers $(d(s_i,s_j))_{i,j}$ are all computable,
uniformly in $i,j$.
\end{definition}

Symbolic spaces, euclidean spaces, functions spaces and manifolds with a
suitable metrics can be endowed with the structure of computable metric
spaces. See for example \cite{G93, HR07, GHR07}.

If $(X,d,\S )$ and $(X^{\prime},d^{\prime},\S ^{\prime})$ are two computable
metric spaces, then the product $(X\times X^{\prime},d_\times,\S \times \S %
^{\prime})$ with $d_\times((x,x^{\prime}),(y,y^{\prime}))=\max(d(x,y),d^{%
\prime}(x^{\prime},y^{\prime}))$ is a computable metric space.

The numbered set of ideal points $(s_i)_i$ induces the numbered set of \emph{%
\textbf{ideal balls}} $\mathcal{B}:=\{B(s_i,q_j):s_i \in S, q_j \in\mathbbm{Q}%
_{>0}\}$. We denote by $B_\uple{i,j}$ the ideal ball $B(s_i,q_j)$.

Let $(X,d,\S )$ be a computable metric space. The computable structure of $X$
assures that the whole space can be ``reached'' using algorithmic means.
Since ideal points (the finite objects of $\S $) are dense, they can
approximate any $x$ at any finite precision. Then, every point $x$ has a
neighborhood basis consisting of ideal balls, denoted $\mathcal{B}%
(x)=\{B\in \mathcal{B}:x\in B\}$ and called its \emph{\textbf{ideal
neighborhood basis}}.

\begin{definition}[Computable points]
A point $x\in X$ is said to be \emph{\textbf{computable}} if its ideal
neighborhood basis $\mathcal{B}(x)$ is r.e.
\end{definition}

\begin{remark} As in the case of reals we have the following characterization: \label{fasttt} $x$ is computable if and only if there is a (total) recursive function $\varphi $ such that $d(s_{\varphi (i)},x)<2^{-i}$.
\end{remark}

Ideal balls are also useful to describe open sets.

\begin{definition}[Recursively open sets]
We say that the set $U\subset X$ is \emph{\textbf{recursively open}} if
there is some r.e set $A$ of ideal balls such that $U=\bigcup_{B\in A}B$.
That is, if there is some r.e set $E\subseteq \mathbbm{N}$ such that $%
U=\bigcup_{i\in E}B_i$.
\end{definition}

We remark that the collection of r.e. open sets can be algorithmically
enumerated.

\begin{definition}
Let $(U_{n})_{n}$ be a sequence of r.e.~open sets. We say that the sequence
is \emph{\textbf{uniformly r.e.}}~or that $U_{n}$ is r.e.~open \emph{\textbf{%
uniformly in $n$}} if there exists an r.e.~set $E\subset \mathbbm{N}^{2}$
such that for all $n$ we have $U_{n}=\bigcup_{i\in E_{n}}B_{i}$, where $%
E_{n}=\{i:(n,i)\in E\}$.
\end{definition}

\begin{examples}
\item Let $(U_n)_n$ be a sequence of open sets such that $U_n$ is uniformly
recursively open. Then the union $\bigcup_n U_n$ is a recursively open set.

\item The universal recursive function $\varphi_u$ induces an enumeration of
the collection $\mathcal{U}$ of all the recursively open sets. Indeed,
define $E:=\{(e,\varphi_u({\langle e,n \rangle})): e,n \in \mathbbm{N}\}$.
Then $\mathcal{U}=\{U_e:e\in N\}$ where $U_e=\bigcup_{i\in E_e}B_i$.

\item The numbered set $\mathcal{U}$ is closed under finite unions and
finite intersections. Furthermore, these operations are \emph{effective} in
the following sense: there exists recursive functions $\varphi^{\cup}$ and $%
\varphi^{\cap}$ such that for all $e,e^{\prime}\in \mathbbm{N}$, $U_e\cup
U_{e^{\prime}}=U_{\varphi^{\cup}({\langle e,e^{\prime}\rangle})}$ and the
same holds for $\varphi^{\cap}$. Equivalently: $U_e\cup U_{e^{\prime}}$ is
recursively open uniformly in ${\langle e,e^{\prime}\rangle}$ (see~\cite{HR07} e.g.).
\end{examples}

\begin{definition}[Computable functions]
\label{functions} A function $T: X \rightarrow Y$ is said to be \emph{%
\textbf{computable}} if $T^{-1}(U^Y_e)$ is recursively open uniformly in $e$.
\end{definition}

It follows that computable functions are continuous. Since we will work with
functions which are not necessarily continuous everywhere, we shall consider
functions which are computable on some subset of $X$. More precisely:
 \begin{definition} A function $T$ is said to be \emph{\textbf{computable on C}} ($C\subset X$) if there is  $U_{n}^{X}$
recursively open uniformly in $n$ such that
 $$T^{-1}(B_{n}^{Y})\cap C=U_{n}^{X}\cap C.$$ The set $C$ is called the \emph{\textbf{
domain of computability}} of $T$.
\end{definition}

As an example we show that a monotone real function whose values over the
rationals are computable, is computable everywhere. This Lemma will also be
used later.

\begin{lemma}
\label{upup}If $f:[0,1]\rightarrow \lbrack 0,1]$ is increasing and $f(r)$
can be computed uniformly, for each rational $r$ then $f$ is computable.
\end{lemma}

\begin{proof}
Let $a,q\in \mathbbm{Q}$. We remark that $f^{-1}((p,q))=\cup _{f(a)\geq
p,f(b)\leq q}(a,b)$ this allows to find a r.e. cover of the interval $%
f^{-1}((p,q))$. The case of a general r.e. open set is straightforward.
\end{proof}

\begin{definition}[Lower semi-computable functions]
\label{lsfunctions} A function $f:X\rightarrow \overline{\mathbbm{R}}$ is said to be \emph{\textbf{lower semi-computable}} if $f^{-1}(q_{n},\infty )$ is recursively open uniformly in $n$.
\end{definition}

It is known that there exists a recursive enumeration of all lower
semi-computable functions $(f_i)_i\geq 0$. 
From the definition follows that lower semi-computable functions are lower semi-continuous. 
\emph{\textbf{Lower semi-computability on D}} is defined as for computable
functions. A function $f$ is \emph{\textbf{upper semi-computable}} if $-f$
is lower semi-computable. It is easy to see that a real function $f$ is
computable if and only if it is upper and lower semi-computable.

Given a probability measure $\mu$, we say that a function is \emph{\textbf{(lower semi-) computable almost everywhere}} if its domain of computability has $\mu$-measure one.


\subsection{Recursively compact sets: approximation from above}\label{s.rec.compact}
We will give some general results about solutions of equations concerning functions computable on some subset. As in many other mathematical situations, to prove the existence of certain solutions we are helped by a suitable notion of compactness. In order to the solution be computable, we will need a recursive version of compacity. Roughly, a compact set is recursively compact if the fact that it is covered by a finite collection of ideal balls can be tested algorithmically (for equivalence with the $\epsilon$-net approach see definition \ref{pree} and proposition \ref{pree2} ).
This kind of notion and the related basic results are already present in the literature in various forms, or particular cases, we give a very compact self contained introduction based on the previously introduced notions.

\begin{definition}
A set $K\subseteq X$ is \defin{recursively compact} if it is compact and there is a recursive
function $\varphi :\mathbbm{N}\rightarrow \mathbbm{N}$ such that $\varphi ({\langle i_{1},\ldots ,i_{p}\rangle })$ halts if and only if $(B_{i_{1}},\ldots ,B_{i_{p}})$ is a covering of $K$.
\end{definition}

\begin{remark}Let $U_{i}$ be the collection of r.e open sets (with its uniform enumeration). It is easy to see that a set $K$ is recursively compact iff $K\subseteq U_{i}$ is
semi-decidable, uniformly in $i$.
\end{remark}

Here are some basic properties of recursively compact sets:

\begin{proposition}\label{p.compact-basic}Let $K$ be a recursively compact subset of $X$.

\begin{enumerate}
\item A singleton $\{x\}$ is recursively compact if and only if $x$ is a computable point.
\item\label{core_compact} If $K'$ is rec. compact then so is $K\cup K'$.
\item if $U$ is recursively open, then $K'=K\setminus U$ is rec compact.
\item The diameter of  $K$ is upper semi-computable.
\item The distance to $K$ : $d_{K}(x):=\inf \{d(x,y): y \in K  \}$ is lower-computable
\item\label{supinf_compact} If $f:X\to \R$ is lower-computable then so is $\inf_{K}f$ 
\item if $f:X\to \R$ is upper-computable then so is $\sup_{K}f$
\end{enumerate}
\end{proposition}

\begin{proof}(1) A point $x$ is computable iff $x\in U_{i}$ is semi-decidable uniformly in $i$. (2) $K\cup K' \subset U$ iff $K \subset U$ and $K' \subset U$. (3) Remark that $K\setminus U\subseteq V\iff K\subseteq U\cup V$ and $U\cup V$ is recursively open uniformly in $U$ and $V$. (4) $\mbox{diam} K=\inf \{q: \exists s, K\subseteq B(s,q)\}$. (5) For $x\in X$ and $q\in \Q$ define $U_{q,x}:=\{y : d(x,y)>q\}$, which is a constructive (in $x$) open set.  Then $d_{K}(x)=\sup\{q : K \subset U_{q,x}\}$ is lower-computable. (6) $\inf_K f=\sup\{q:K\subseteq f^{-1}(q,+\infty)\}$. (7) $\sup_K f=\inf\{q:K\subseteq f^{-1}(-\infty,q)\}$.




%

%



%

\end{proof}

\begin{remarks}
\item The arguments are uniform. In point 1) for instance, this means that there is an algorithm which takes a program computing $x$ and outputs a program testifying the rec. compacity of $\{x\}$, and vice-versa.
\item When $X$ itself is rec. compact, a subset $K$ is rec. compact iff $d_{K}$ is lower-computable. Indeed, $K=X\setminus \{x : d_{K}(x)>0\}$.
\end{remarks}

\begin{corollary}\label{c.countable-intersection} If $(K_{i})_{i\in \N}$ are uniformly recursively compact sets, then so is $\bigcap_{i\in \N} K_{i}$.
\end{corollary}
\begin{proof}The complements of recursively compact sets are r.e open. Then by proposition \ref{p.compact-basic}, part (\ref{core_compact}) the set $\bigcap_{i\in \N} K_{i}=K_{0}\setminus (\bigcup_{i>0}K_{i}^c)$ is recursively compact.
\end{proof}

It is important to remark that a recursively compact set needs not contain computable points. This will be used in section \ref{lst}.
\begin{proposition}\label{p.no-comp}
There exists a nonempty recursively compact set $K\subset [0,1]$ containing no computable points. 
\end{proposition}
\begin{proof}
Let $I_{n}$ be an enumeration of all the rational intervals and $\epsilon>0$ be a rational number. Put $E=\{i\geq 1:\varphi_i(i)\mbox{ halts and }|I_{\varphi_i(i)}|<\epsilon 2^{-i}\}$. $E$ is a r.e. subset of $\N$. Let $U=\bigcup_{i\in E}I_i$: $\lambda(U)\leq \sum_{i\in E}\epsilon 2^{-i}\leq \epsilon$. Let $x\in[0,1]$ be a computable real number. There is a total recursive function $\varphi_i$ such that $|I_{\varphi_i(n)}|<\epsilon 2^{-n}$ and $x\in I_{\varphi_i(n)}$ for all $n$. Hence $i\in E$, so $x\in U$. Hence $U$ contains all computable points. As $[0,1]$ is recursively compact, so is $K=[0,1]\setminus U$.
\end{proof}

Now we start to show that many statements about topology and calculus on
metric spaces can be easily translated to the computable setting: the first one
says that the image of a recursively compact is still a recursively compact.

\begin{proposition}[Stability by computable functions]\label{p.stability}
Let $f:K\subseteq X\to Y$ be a computable function defined on a recursively
compact set $K$. Then $f(K)$ is recursively compact.
\end{proposition}

\begin{proof}
Indeed, $f(K)\subseteq U\iff K\subseteq f^{-1}(U)$. As $f^{-1}(U_e)\cap
K=U_{\varphi(e)}\cap K$ where $\varphi$ is a total recursive function, $f(K)\subseteq U_e\iff K\subseteq U_{\varphi(e)}$.
\end{proof}

Remark that the argument is uniform: if $(K_i)_{i\in\mathbbm{N}}$ is a
sequence of uniformly recursively compact subsets of $X$ on which $f$ is
defined, then $(f(K_i))_{i\in\mathbbm{N}}$ is a sequence of uniformly
recursively compact subsets of $Y$. We will say that $f(K)$ is recursively
compact \emph{uniformly in $K$}.

As a first simple example of application, we observe that in some cases the global attractor of a (computable) dynamical system can be approximated by an algorithm to any given accuracy.

\begin{corollary}\label{c.rec.compact.attractor}
Let $X$ be a recursively compact computable metric space and $T$ a computable dynamics on it. Then the set:
$$
\Lambda:=\bigcap_{n\geq 0}T^n (X)
$$
is recursively compact.
\end{corollary}
\begin{proof}By proposition \ref{p.stability} and corollary \ref{c.countable-intersection}
\end{proof}

We remark that these and other frameworks of  ``exact computability and rigorous
approximation'' have been previously used to study the computability of several similar objects such as Julia or Mandelbrot sets (\cite{H05, BY06, BBY06, BBY07}, \cite{Del97}), or the existence and some basic properties of Lorentz
attractor (\cite{Tuc99}).

Here is a computable version of Heine's theorem.

\begin{definition}
A function $f:X\rightarrow Y$ \ between metric spaces is \defin{recursively
uniformly continuous} if there is a recursive $\delta :\mathbbm{Q\rightarrow Q}
$ such that for all $\epsilon>0$, $\delta(\epsilon)>0$ and $\forall x\in X$,
\begin{equation}
f(B(x,\delta (\epsilon )))\subset B(f(x),\epsilon ).
\end{equation}
\end{definition}

\begin{proposition}
Let $X$ and $Y$ be two computable metric spaces. Let $K\subseteq X$ be
recursively compact and $f:K\rightarrow Y$ be a computable function. Then $f$
is recursively uniformly continuous.
\end{proposition}

\begin{proof}
First, $K\times K$ is a recursively compact subset of $X\times X$. For each
rational number $\epsilon >0$, define $U(\epsilon )=\{(x,x^{\prime })\in
K^{2}:d(f(x),f(x^{\prime }))<\epsilon \}$ and $K(\epsilon )=K\times
K\setminus U(\epsilon )$: they are respectively recursively open and
recursively compact, uniformly in $\epsilon $. Hence, the function $\delta
(\epsilon ):=\inf \{d(x,y):(x,y)\in K(\epsilon )\}$ is lower semi-computable
(proposition \ref{supinf_compact}).

Now, $f$ is uniformly continuous if and only if $\delta(\epsilon)>0$ for
each $\epsilon>0$. By the classical Heine's theorem, this is the case, so by
lower semi-computability of $\delta(\epsilon)$, one can compute from $%
\epsilon$ some positive $\delta\leq \delta(\epsilon)$. %
%
%
%
\end{proof}

\begin{theorem}
\label{computable_zero} Let $K$ be a recursively compact subset of $X$ and $f:K\to \mathbbm{R}$ be a computable function. Then every isolated zero of $f$
is computable.
\end{theorem}

\begin{proof}
Let $x_0$ be an isolated zero of $f$. Let $s,r$ be an ideal point and a
positive rational number such that $x\in B(s,r)$ and the only zero of $f$
lying in $\overline{B}(s,r)$ is $x_0$. The set $N=\{x:f(x)\neq 0\}\cup
\{x:d(x,s)>r\}$ is recursively open in $K$ (that is, $N\cap K=U\cap K$ with $U$ recursively open), so $\{x_0\}=K\setminus N=K\setminus U$ is
recursively compact by proposition \ref{p.compact-basic}. Hence, $x_0$ is a
computable point.
\end{proof}

\begin{remark}\label{rmk3}Observe that the argument is uniform in $f$ and an ideal ball isolating the zero. In particular, there is an algorithm which takes a finite description of $f$ and the ball   and outputs his zero if it is unique.
\end{remark}

\begin{corollary}
\label{computable_fixed} Let $K$ be a recursively compact subset of $X$ and $f:K\to X$ be a computable function. Then every isolated fixed point of $f$
is computable.
\end{corollary}

\begin{proof}
Apply the preceding theorem to the function $g:X\to\mathbbm{R}$ defined by $g(x)=d(x,f(x))$.
\end{proof}


\subsection{Recursively precompact}
In this subsection we prove the equivalence between the notion of recursive compactness given above and another natural approach (which will be used later) to recursive compactness, where it is supposed the existence of an algorithm to construct   $\epsilon $-nets. 

\begin{definition}\label{pree}
A CMS is \defin{recursively precompact} if there is a total recursive function $\varphi:\mathbbm{N}\to\mathbbm{N}$ such that for all $n$, $\varphi(n)$
computes a $2^{-n}$-net: that is $\varphi(n)={\langle i_1,\ldots,i_p \rangle}
$ where $(s_{i_1},\ldots,s_{i_p})$ is a $2^{-n}$-net.
\end{definition}

Here is a computable version of a classical theorem:

\begin{proposition}\label{pree2}
Let $X$ be a CMS. $X$ is recursively compact if and only if it is complete
and recursively precompact.
\end{proposition}

\begin{proof}
If $X$ is recursively compact then we define the following algorithm: it
takes $n$ as input, then enumerates all the ${\langle i_1,\ldots,i_p \rangle}
$, and tests if $(B(s_{i_1},2^{-n}),\ldots,B(s_{i_p},2^{-n}))$ is a covering
of $X$ (this is possible by recursive compacity). As $X$ is compact, hence
precompact, such a covering exists and will be eventually enumerated: output
it. The algorithm makes $X$ recursively precompact.

Suppose that $X$ is complete and recursively precompact. Let $(B(s_1,q_1),\ldots,B(s_k,q_k))$ be ideal balls: we claim that $(B(s_1,q_1),\ldots,B(s_k,q_k))$ covers $X$ if and only if there exists $n$
such that each point $s$ of the $2^{-n}$-net given by recursive
precompactness lies in a ball $B(s_i,q_i)$ satisfying $d(s,s_i)+2^{-n}<q_i$.
The procedure which enumerates all the $n$ and semi-decides this halts if
and only if the initial sequence of balls covers $X$. We leave the proof of
the claim to the reader (take $n$ such that $2^{-n}$ is less than the
Lebesgue number of the finite covering).
\end{proof}

The following observation is also worth noticing.

\begin{proposition}
Let $X$ be a computable metric space. If $X$ (as a subset of $X$) is
recursively compact, then the set $C(X)$ of continuous functions from $X$ to 
$\mathbbm{R}$ with the distance induced by the uniform norm is a computable
metric space.

The function $eval:\mathcal{C}(X)\times X\to\mathbbm{R}$ mapping $(f,x)$ to $f(x)$ is computable.

Let $Y$ be a computable metric space: for every computable function $f:Y\times X\to \mathbbm{R}$, the function $Y\to C(X)$ mapping $y$ to $f_y:x\mapsto f(y,x)$ is computable.
\end{proposition}

\subsection{Recursively closed sets: approximable from below}

From the computability viewpoint, the properties of recursively closed sets are, in a sense, complementary to those of recursively compact sets.

\begin{definition}
A closed set $F$ is \emph{\textbf{recursively closed}} if the set $\{B(s,r):
B(s,r)\cap F\neq \emptyset\}$ is r.e.
\end{definition}

A closed set $F$ is recursively closed if $F\cap U$ is semi-decidable for
r.e open sets $U$. It is easy to see that the union of two recursively closed
sets is also recursively closed. The closure of any recursively open set is
recursively closed: $B\cap \overline{U}\neq \emptyset\iff \exists s\in B\cap
U$.

The following proposition will be used later.

\begin{proposition}\label{p.comp-closed}
Let $F$ be a recursively closed subset of $X$. Then there exists a sequence of uniformly computable  points $x_{i}\in F$ which is dense in $F$.
 \end{proposition}
\begin{proof}

Since $\{n\in \N : B_n=B(s_{n},q_{n})\cap F \neq \emptyset\}$ is r.e, given some ideal ball $B=B(s,q)$ intersecting $F$, the set $\{n \in \N : \overline{B_n}\subset B, q_{n}\leq 2^{-n}, B_n\cap F \neq \emptyset \}$ is also r.e. Then we can effectively construct an exponentially decreasing sequence of ideal balls intersecting $F$. Hence $\{x\}=\cap_k B_k$ is a computable point lying in $F$.

\end{proof}




We remark that by this, Proposition \ref{p.no-comp} shows a recursive compact which is not recursively closed.
For the sake of completeness, let us state some useful simple properties.

\begin{proposition}Let $F$ be a recursively closed subset of $X$. Then:
\begin{itemize}
\item[(1)] The diameter of $F$ is lower semi-computable,
uniformly in $F$.
\item[(2)] If $f:F\to \mathbbm{R}$ is lower semi-computable, then so is $\sup_F f$. 
\item[(3)] If $f:F\to \mathbbm{R}$ is upper semi-computable, then so is $\inf_F f$. 
\end{itemize}
\end{proposition}

\begin{proof}
(1) Let $C(s,q)$ be the complement of the closed ball $\overline{B}(s,q)$, that
is $C(s,r)=\{x:d(x,s)>q\}$: this is a recursively open set, uniformly in $%
s,q $. Then $\mbox{diam} F=\sup \{q:\exists s,C(s,r)\cap F\neq \emptyset\}$. (2)  $\sup_F f=\sup
\{q:f^{-1}(q,+\infty)\cap F\neq \emptyset\}$. (3) Apply (2) to $-f$.
%
%
\end{proof}

\begin{corollary}
Let $K$ be  recursively closed and recursively compact subset of $X$. If $f:K\to \mathbbm{R}^+$ is a computable function, then so are $\inf_K f$ and $\sup_K f$. 
\end{corollary}



\section{Computable measures\label{meas}}

Let us consider the space $PM(X)$ of Borel probability measures over $X$. We
recall that $PM(X)$ can be seen as the dual of the space $C_{0}(X)$ of
continuous functions with compact support over $X$ and recall the notion of
weak convergence of measures:

\begin{definition}
$\mu _{n}$ is said to be \defin{weakly convergent} to $\mu $ if $\integral{f}{\mu_n}\rightarrow \integral{f}{\mu}$ for each $f\in C_{0}(X)$.
\end{definition}

Let us introduce the Wasserstein-Kantorovich distance between measures. Let $\mu _{1}$ and $\mu _{2}$ be two probability measures on $X$ and consider:

\begin{equation*}
W_{1}(\mu _{1},\mu _{2})=\underset{f\in 1\text{-Lip}(X)}{\sup }\left|\integral{f}{\mu_1}-\integral{f}{\mu _2}\right|
\end{equation*}

where $1\mbox{-Lip}(X)$ is the space of 1-Lipschitz functions on $X$. We remark
that since adding a constant to the test function $f$ does not change the
above difference $\integral{f}{\mu _1}-\integral{f}{\mu _2}$ then the supremum can be taken over the set of 1-Lipschitz functions mapping a distinguished ideal point $s_0$ to $0$. 
The distance $W_{1}$ has moreover the
following useful properties which will be used in the following

\begin{proposition}[\cite{AGS} Prop 7.1.5]\label{ambros}\mbox{}

\begin{enumerate}
\item $W_{1}$ is a distance and if $X$ is bounded, separable and complete, then $PM(X) $ with this distance is a separable and complete metric space.

\item If $X$ is bounded, a sequence is convergent for the $W_{1}$ metrics if
and only if it is convergent for the weak topology.

\item If $X$ is compact $PM(X)$ is compact with this topology.
\end{enumerate}
\end{proposition}

Item (1)  has an effective version: $PM(X)$ inherits the computable metric structure of $X$. Indeed, given the set $S_{X}$ of ideal points of $X$ we can naturally define a set
of ideal points $S_{PM(X)}$ in $PM(X)$ by considering finite rational convex
combinations of the Dirac measures $\delta _{x}$ supported on ideal points $x\in S_{X}$. This is a dense subset of $PM(X)$. The proof of the following proposition can be found in (\cite{HR07})

\begin{proposition}
If $X$ bounded then $PM(X)$ with the $W_{1}$ distance (and $S_{PM(X)}$ as a set of ideal points) is a computable metric space.
\end{proposition}

A measure $\mu $ is then computable if there is a fast sequence $(\mu
_{n})\in S_{PM(X)}$ converging to $\mu $ (see remark \ref{fasttt}) in
the $W_{1}$ metric (and hence for the weak convergence).

Now, point (3) of proposition \ref{ambros} also has an effective version:

\begin{lemma}
\label{recomp}If $X$ is a recursively precompact metric space, then $PM(X)$ with the $W_{1}$ distance is a recursively precompact metric space.
\end{lemma}

\begin{proof}
We will show how to effectively find an $r-$net for each $r$ of the form $r=%
\frac{1}{n},n\in \mathbbm{N}$. Let us consider the set $S_{r}=\{\frac{k}{n}%
,0\leq k\leq n\}$ subdividing the unit intervals in equal segments. Let us
also consider an $r$-net $N_{r}=\{x_{1},...x_{m}\}$ constructed by recursive
compactness of $X.$ Now let us consider the set $\Upsilon _{r}$ of measures
with support in $N_{r}$ given by 
\begin{equation*}
\Upsilon _{r}=\{k_{1}\delta _{x_{1}}+...+k_{m}\delta _{x_{m}}~s.t.~k_{i}\in
S_{r}~,~k_{1}+...+k_{m}=1\}.
\end{equation*}%
This is a \thinspace $2r$ net in $PM(X)$. To see this let us consider a
measure $\mu $ on $X$ and a ball $B(x_{1},r)$ centered in $x_{1}\in X$. Let
us consider the measure $\mu _{1}$ defined by

\begin{equation*}
\mu _{1}(A)=\mu (A)-\mu (B(x_{1},r)\cap A)+\delta _{x_{1}}(A)
\end{equation*}
for each measurable set $A\subset X$. The measure $\mu _{1}$ is obtained
transporting the mass contained in the ball $B(x_{1},r)$ to its center. Then 
$W_{1}(\mu _{1},\mu )\leq r\mu (B(x_{1},r))$. Let us now consider the
sequence of measures $\mu _{1},...,\mu _{m}$ where $\mu _{1}$ is as before
and the other ones are given by%
\begin{equation*}
\mu _{i}(A)=\mu _{i-1}(A)-\mu _{i-1}(B(x_{i},r)\cap A)+\delta _{x_{i}}(A),
\end{equation*}%
at the end $\mu _{m}$ is a measure with support in $N_{r}$ and by the
triangle inequality $W_{1}(\mu _{m},\mu )\leq r.$

Now $\mu _{m}$ has the same support as the measures in $\Upsilon _{r}$ and
there is $\nu \in \Upsilon _{r}$ such that $|\integral{f}{\mu_n}-\integral{f}{\nu}|\leq r$ for each $f\in 1\mbox{-Lip}(X)$, hence $W_{1}(\mu _{m},\nu )\leq r$ and
then $W_{1}(\mu ,\nu )\leq 2r$ and this proves the statement.
\end{proof}

We now use the recursive enumeration of lower semi-computable functions $(f_{i})_{i}\geq 0$ to characterize computability on $PM(X)$ (see \cite{HR07} corollary 4.3.1):

\begin{lemma}
\label{comp_meas}Let $X$ be a bounded computable metric space and $\mathcal{S%
}$ be any subset of $PM(X)$, then:

\begin{enumerate}
\item $\mu \in PM(X)$ is computable iff the function $\mu \mapsto
\integral{f_i}{\mu}$ is lower semi-computable, uniformly in $i$,

\item $L:PM(X)\rightarrow PM(X)$ is computable on $\mathcal{S}$ iff the
function $\mu \mapsto \integral{f_i}{L(\mu)}$ is lower semi-computable on $\mathcal{S}$, uniformly in $i$.
\end{enumerate}
\end{lemma}

This gives:

\begin{lemma}
\label{comp_D} If $g_{i}:X\rightarrow \mathbbm{R}^{+}$ is a uniform sequence
of functions which are lower semi-computable on $X\setminus D$, then $\mu
\mapsto \integral{g_i}{\mu}$ is lower semi-computable on 
\begin{equation}
PM_{D}(X):=\{\mu :\mu (D)=0\}
\end{equation}
uniformly in $i$.
\end{lemma}

\begin{proof}
For each $i$, one can construct a lower semi-computable function $\hat{g}_i$ satisfying $\hat{g}_{i}=g_{i}$ on $X\setminus  D$ (see \cite{HR07}, subsection 3.1). Since the function $\mu \mapsto \integral{\hat{g}_i}{\mu}$
is lower semi-computable, uniformly in $i$ and $\mu (D)=0$, we have that on $PM_{D}(X) $ it coincides with $\mu \mapsto \integral{\hat{g}_i}{\mu}$, which is
then lower semi-computable on $PM_{D}(X)$, uniformly in $i$.
\end{proof}

An interesting remark about computable measures is that they must have computable points in the support. This will be used in section \ref{example2}. 
\begin{proposition}\label{compoint}
If $\mu$ is a computable probability measure, then there exists computable points in the support of $\mu$.
\end{proposition}
\begin{proof}The sequence of functions $f_{i}:=1_{B_{i}}$ (the indicator functions of ideal balls) are uniformly lower semi-computable. By lemma \ref{comp_meas}, the numbers $\integral{f_{i}}{\mu}=\mu(B_{i})$ are uniformly lower semi-computable.  Hence, the set $\{B_{i}: \mu(B_{i})>0\}$ is recursively enumerable. In other words, the support of $\mu$ is a recursively closed set. Proposition \ref{p.comp-closed} allows to conclude.
\end{proof}


\section{Dynamical systems, statistical behavior, invariant measures}\label{DS}

Let $X$ be a metric space, let $T:X\mapsto X$ be a Borel measurable map. Let $\mu $ be
an invariant measure. A set a $A$ is called $T$-invariant if $%
T^{-1}(A)=A(mod0)$. The system $(T,\mu )$ is said to be ergodic if each $T$%
-invariant set has total or null measure. In such systems the famous
Birkhoff ergodic theorem says that time averages computed along $\mu $
typical orbits coincides with space average with respect to $\mu .$ More
precisely, for any $f\in L^{1}(X)$ it holds 
\begin{equation}
\underset{n\rightarrow \infty }{\lim }\frac{S_{n}^{f}(x)}{n}=\int fd\mu ,
\label{Birkhoff}
\end{equation}%
for $\mu $ almost each $x$, where $S_{n}^{f}=f+f\circ T+\ldots +f\circ
T^{n-1}.$

This shows that in an ergodic system, the statistical behavior  of observables, under typical realizations of the system is given by the average of the observable made with the invariant measure.

In case $X$ is a manifold (possibly with boundary). We say that a point $x$ belong to the basin of an invariant measure $\mu $ if Equation \ref{Birkhoff} holds at $x$ for each continuous $f$ (the average on the $x$ orbit represent the average under the measure). An SRB measure is an invariant measure having a positive Lebesgue measure basin (for more details and a general survey see \cite{Y}).

In the applied literature the most common method to simulate or understand the above statistical behaviors is to compute and study some trajectory. This method has three main theoretical problems which motivates the search of another approach: 
\begin{itemize}
\item numerical error, 
\item tipicality of the sample,
\item how many sample points are necessary?
\end{itemize}
the first (and widely known) problem is the amplification of the
numerical error (if the system is sensitive to initial conditions as most
interesting systems are). Here the shadowing results are often invoked to
justify the correctness of simulations, but rigorous results are proved only
for a small class of systems (see e.g.\cite{Pal00}) and moreover the mere
existence of a shadowing orbit does not say anything about its typicality
(see e.g. \cite{Bla89, Bla94} for a further discussion on numerical errors).

The second problem is indeed that this method should compute, in order to be
useful, a trajectory which shows the ``typical'' behavior of the system: a
behavior which take place with large or full probability in some sense. The main
problem here is the fact that the set of initial conditions the computer has
access to, being countable, has probability zero. Hence, there is no
guarantee that what we see on the screen is typical in some sense.
On the contrary, in a chaotic system, typical orbits are far from being describable by a finite program. It is true for example that in an ergodic system having positive entropy $h$ a  typical $n$ step orbit segment needs approximatively a program which is $hn$ bits long to be described (up some approximation $\epsilon $, see e.g. \cite{B83} for the original result or \cite{Ga00} and \cite{GHR09c} for a version in the framework of computable analysis). We remark, however, that if one looks for points which  behave as typical for Birkhoff averages (hence they behave as typical for some given particular aspect) there are some rigorous results partly supporting this way to proceed: in several classes of systems there are computable initial conditions which behave as typical with respect to Birkhoff averages (see \cite{GHR07} for a precise result). 

The third problem however remains. Even if you find a program describing a typical orbit of the system: how many iterations should be considered to be near to the limit behavior, so that the orbit represents the invariant measure up to a certain approximation? although this problem can be approached rigorously in some cases (see \cite{CCS} e.g.) we will not adopt this point of view.
 We will study the system's statistical behavior by directly computing the invariant measure as fixed points of a certain transfer operator.

\subsection{The transfer operator\label{PF}}

A function $T$ between metric spaces naturally induces a function between probability measure spaces. This function $L_T$ is linear and is called transfer operator (associated to $T$). Measures which are invariant for $T$ Invariant  measures are fixed points of  $L_T$. 
 
Let us consider a computable metric space $X$ endowed with a Borel
probability measure $\mu$ and with a dynamics defined by a
measure-preserving function $T:X\rightarrow X$. Let us also consider the
space $PM(X)$ of Borel probability measures on $X.$

Let us define the function $L_{T}:PM(X)\rightarrow PM(X)$ by duality in the
following way: if $\mu \in PM(X)$ then $L_{T}(\mu )$ is such that

\begin{equation*}
\integral{f}{L_{T}(\mu )}=\integral{f\circ T}{\mu}
\end{equation*} for each $f\in C_{0}(X)$.
In next sections, invariant measures will be found as solutions of the
equation $W_{1}(\mu ,L(\mu ))=0.$ To apply Theorem \ref{computable_zero}
and Corollary \ref{computable_fixed} to this equation we need that $L$ is computable.
We remark that if $T$ is not continuous then $L$ is not necessarily
continuous (this can be realized by applying $L$ to some delta measure
placed near a discontinuity point) hence not computable. Still, we have
that $L$ is continuous (and its modulus of continuity is computable) at all
measures $\mu $ which are \textquotedblleft far enough\textquotedblright\
from the discontinuity set $D$. This is technically expressed by the
condition $\mu (D)=0$.

We remark that with the general tools introduced before, the proof is immediate.

\begin{theorem}
\label{L_comp} Let $X$ be a computable metric space and $T:X\rightarrow X$
be a function which is computable on $X\setminus D$. Then $L_{T}$ is
computable on the set of measures 
\begin{equation}
PM_{D}(X):=\{\mu \in PM(X):\mu (D)=0\}.
\end{equation}
\end{theorem}

\begin{proof}
Note that if $f$ is lower semi-computable, then $f\circ T$ is lower
semi-computable on $X\setminus D$. The result then follows from lemmas \ref{comp_D} and %
\ref{comp_meas}.
\end{proof}

In particular, if $T$ is computable on the whole space $X$ then $L$ is
computable on all $PM(X)$.


\subsection{Computing invariant ``regular'' measures}

The above tools allow to ensure the computability of $L_{T}$ on a large
class of measures. This will allow to apply Corollary \ref{computable_fixed}
and see an invariant measure as a fixed point.

\begin{theorem}
\label{comp_inv} Let $X$ be a computable metric space and $T$ be a function
which is computable on $X\setminus D$. Suppose there is a recursively
compact set of probability measures $V\subset M(X)$ such that for every $\mu
\in V$, $\mu (D)=0$ holds. Then every invariant measure isolated in $V$ is
computable.
\end{theorem}

\begin{proof}
By Theorem \ref{L_comp}, $L_{T}$ is computable on $V$. Since $V$ is
recursively compact, theorem \ref{computable_zero} allows to compute any
invariant isolated measure in $V$ as a solution of the equation $L_T(\mu)=\mu $.
\end{proof}

\begin{remark}\label{comprem}
This theorem is uniform: there is an algorithm which takes as inputs finite descriptions of $T, V$ and an ideal ball in $M(X)$ which isolates an invariant measure $\mu$, and outputs a finite description of $\mu$ (see the above proof and Remark \ref{rmk3}). 
\end{remark}
A trivial consequence of  Theorem \ref{comp_inv} is the following:

\begin{corollary}
If a computable system  as above is uniquely ergodic and its invariant measure $\mu$ satisfy $\mu (D)=0$, then it is a computable measure.
\end{corollary}

The main problem in the application of theorem \ref{comp_inv} is the
requirement that the invariant measure we are trying to compute, should be
isolated in $V$. In general the space of invariant measures in a given
dynamical system could be very large (an infinite dimensional convex in $%
PM(X)$ ) to isolate a particular measure we can restrict and consider a
subclass of "regular" measures.

Let us consider the following \emph{seminorm}:

\begin{equation*}
\left\Vert \mu \right\Vert _{\alpha }=\sup_{x\in X,r>0}\frac{\mu (B(x,r))}{%
r^{\alpha }}.
\end{equation*}

\begin{proposition}
\label{1}If $X$ is recursively compact then 
\begin{equation}
V_{\alpha ,K}=\{\mu \in PM(X):\left\Vert \mu \right\Vert _{\alpha }\leq K\}
\end{equation}%
is recursively compact.
\end{proposition}

\begin{proof}$U=\{\mu \in PM(X):\left\Vert \mu \right\Vert
_{\alpha }>K\}$ is recursively open. Indeed, $\norm{\mu}_\alpha > K$ iff there exists $s,r\in \S\times \Q$ for which $\mu(B(s,r))>qr^{\alpha}$. As $\mu \mapsto \mu(B(s,r))$ is lower semi-computable uniformly in $s,r$, the sets $U_{s,r}:=\{\mu: \mu(B(s,r))>Kr^{\alpha}\}$ are uniformly recursively open subsets of $PM(X)$. Hence, $U=\cup_{s,r}U_{s,r}$ is recursively open.

Now, $V_{\alpha
,K}=PM(X)\setminus U$. As $PM(X)$ is recursively compact (see Lemma \ref{recomp}) and $U$ is recursively
open,
then proposition \ref{p.compact-basic} part (3) allows to conclude.
\end{proof}

In theorem \ref{comp_inv} we require that $\mu (D)=0$ holds. This is
automatically true in many examples when the measure is regular and the set $%
D$ is reasonably small.

\begin{proposition}
\label{2}Let $X$ be recursively compact and $T$ be computable on $X\setminus
D$, with $dim_{H}(D)<\infty $. Then any invariant measure isolated in $%
V_{\alpha ,K}$ with $\alpha >dim_{H}(D)$ is computable.
\end{proposition}

\begin{proof}
Let us first prove that $\mu (D)=0$ for all $\mu \in V_{\alpha ,K}$. For all 
$\epsilon >0$, there is a covering $(B(x_{i},r_{i}))_{i}$ of $D$ satisfying $\sum_{i}r_{i}^{\alpha }<\epsilon $. Hence $\mu (D)\leq \sum_{i}\mu
B(x_{i},r_{i})\leq 2^{\alpha }K\sum_{i}r_{i}^{\alpha }\leq 2^{\alpha
}K\epsilon $. As this is true for each $\epsilon >0$, $\mu (D)=0$.

The result then follows from the fact that $V_{\alpha ,K}$ is recursively
compact and Theorem \ref{comp_inv}.
\end{proof}
\begin{remark}
Once again, this is uniform in $T,\alpha,K$. 
\end{remark}

The above general  proposition allows to obtain as a corollary the computability of many absolutely
continuous invariant measures. For the sake of simplicity, let us consider
maps on the interval.

\begin{proposition}
If $X=[0,1]$, $T$ is computable on $X\setminus D$, with $dim_{H}(D)<1$
and $(X,T)$ has an unique a.c.i.m. $\mu $ with bounded density, then $\mu $
is computable.
\end{proposition}

\begin{proof}

The result follows from the above proposition \ref{2} and the fact that if  $\mu $ is absolutely  continuous and the density of $\mu $ is $f\in L_{1}[0,1]$ then $\left\Vert \mu \right\Vert _{1}=ess sup (f).$ We have to check that there could be not other measures having a finite $1$ norm and not being absolutely continuous.

If we suppose that $\left\Vert \mu \right\Vert_{1}=l$ is finite, then $\mu$ is absolutely
continuous, with bounded density $f\leq l$. Indeed, let us consider the
conditional expectation $E[\mu|I^{n}]$ of $\mu$ to the dyadic $n$-th
grid $I^{n}=\{[k2^{-n},(k+1)2^{-n}),0\leq k\leq2^{n}\}$. 

If $\left\Vert \mu\right\Vert_{1}=l$ a fortiori
implies $0\leq E[\mu|I^{n}]\leq l$ a.e..
By the first Doobs martingale convergence it follows
that $E[\mu|I^{n}]$ has an a.e. pointwise limit $f$ and $f\leq l$ a.e.. Since $f$ is bounded then it is a density for $\mu$.

\end{proof}

$d$-dimensional submanifolds of $\mathbbm{R}^n$  can naturally be endowed with a natural structure of computable metric spaces ( see \cite{GHR07}).
Considering a dyadic grid on $\mathbbm{R}^d$ and chart diffeomorphisms it is straightforward to prove, in the same way as before

\begin{corollary}
Let $X$ be a recursively compact $d$ dimensional $C^{1}$ submanifold of $%
\mathbbm{R}^{n}$ (with or without boundary). If $T$ is computable on $X\setminus D$, with $dim_{H}(D)<d$ and $(X,T)$ has a unique a.c.i.m. $\mu $ with bounded
density, then $\mu$ is computable.
\end{corollary}

As it is well known,  interesting examples of systems having an unique a.c.i.m. (with bounded density as required) are topologically transitive {\em piecewise expanding maps} on the interval or {\em expanding maps} on manifolds (see \cite{V97} for precise definitions). Provided that the dynamics is computable we then have by the above propositions that the a.c.i.m. is computable too.


\subsection{Unbounded densities}

The above results ensure computability of measures having an a.c.i.m. with
bounded density. If we are interested in situations where the density is
unbounded, we can consider a new norm, ``killing'' singularities.

Let us hence consider a computable function $f:X\rightarrow \mathbbm{R}$ \ and

\begin{equation*}
\left\Vert \mu \right\Vert _{f,\alpha }=\sup_{x\in X,r>0}\frac{f(x)\mu
(B(x,r))}{r^{\alpha }}.
\end{equation*}
Propositions \ref{1} and \ref{2} also hold for this norm. If $f$ is such
that $f(x)=0$ when $\lim_{r\rightarrow 0}\frac{\mu (B(x,r))}{r^{\alpha }}%
=\infty $ this can let the norm to be finite when the density diverges.

As an example, where this can be applied, let us consider the Manneville Pomeau maps on the unit interval. These are maps of the type $x\rightarrow x+x^{z}(\func{mod}~1)
$. When $1<z<2$ the dynamics has an unique a.c.i.m. $\mu _{z}$ having density $%
e_{z}(x)$ which diverges in the origin as $e_{z}(x)\asymp x^{-z+1}$ and it
is bounded elsewhere (see \cite{I03} section 10 and \cite{V97} section 3
e.g.). If we consider the norm $\left\Vert .\right\Vert _{f,1}$ with $%
f(x)=x^{2}$ we have that $\left\Vert \mu _{z}\right\Vert _{f,1}$ is
finite for each such $z$. By this it follows that the measure $\mu _{z}$ is
computable.


\section{Computable systems having not computable invariant measures\label{lst}}

We have seen that the technique presented above proves the computability of many a.c.i.m. which are also SRB measures. As we have seen in the introduction, with other techniques it is possible to prove the computability of other SRB measures (axiom A systems e.g., see \cite{GHR07}).
This raises naturally the following question:  a computable systems does necessarily have a computable invariant measure? 
what about ergodic SRB measures?
 
The following is an easy example showing that this
is not true in general even in quite regular systems, hence the whole question of computing invariant measures has some subtlety.

Let us consider a system on the unit interval given as follows. Let $\tau\in (0,1)$ be a lower semi-computable real number which is not computable. There is a computable sequence of rational numbers $\tau_i$ such that $\sup_i \tau_i=\tau$. For each $i$, define $T_i(x)=\max(x,\tau_i)$ and $T(x)=\sum_{i\geq 1} 2^{-i}T_i$. The functions $T_i$ are uniformly computable so $T$ is also computable.

\begin{figure}[h]
\includegraphics[height=5cm]{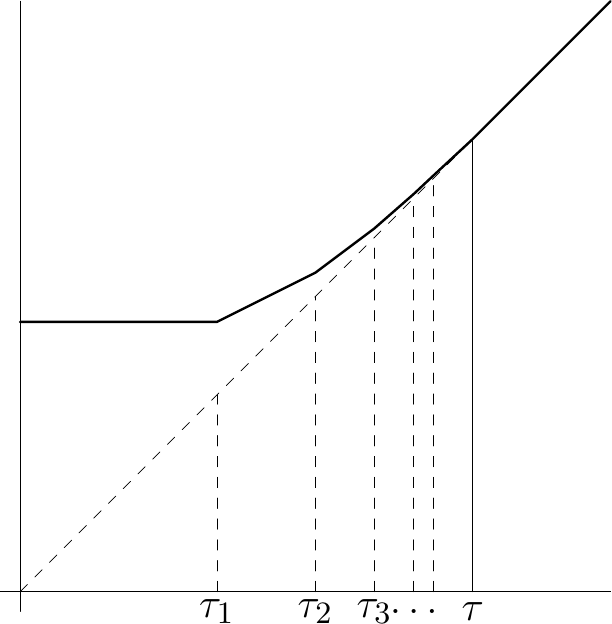}
\caption{The map $T$.}
\end{figure}

Now, $T$ is non-decreasing, and $T(x)>x$ if and only if $x<\tau$.

The system $([0,1],T)$ is hence a computable dynamical
system. This system has a SRB ergodic invariant measure which is $\delta
_{\tau }$, the Dirac measure placed on $\tau $. The measure is SRB because $\tau$ attracts all the interval at its left. Since $\tau $ is not
computable then $\delta _{\tau }$ is not computable. We remark that
coherently with the previous theorems $\delta _{\tau }$ is not isolated.

We remark that by a simple dichotomy argument we can prove that a computable function from $[0,1]$ to itself must have a computable fixed point. Hence it is not possible to construct a system over the interval having no computable invariant measure (we always have the $\delta$ over the fixed point). With some more work we will see that such an example can be constructed on the circle.

\subsection{A computable system having no computable invariant measure}\label{example2}
We go further and exhibit a computable dynamical system on a compact space which has \emph{no} computable invariant probability measure.

We consider the unit circle $S$, identified with $\R/\Z$. It naturally has a computable metric structure inherited from that of $\R$.

On $S$, there is a computable map with no computable invariant probability measure.
We construct such a map $T:[0,1]\to\R$ satisfying $T(1)=T(0)+1$, and consider its quotient on the unit circle.


From proposition \ref{p.no-comp} we know that there is a non-empty recursively compact set $K$ containing no computable point. Let $U=(0,1)\setminus K$: this is a r.e. open set, so there are computable sequences $a_i,b_i$ ($i\geq 1$) such that $0<a_i<b_i<1$ and $U=\bigcup_i (a_i,b_i)$. Let us define non-decreasing, uniformly computable functions $f_i:[0,1]\to[0,1]$ such that $f_i(x)>x$ if $x\in (a_i,b_i)$ and $f_i(x)=x$ otherwise. For instance, $f_i(x)=2x-a_i$ on $[a_i,\frac{a_i+b_i}{2}]$ and $f_i(x)=b_i$ on $[\frac{a_i+b_i}{2},b_i]$.

\begin{figure}[h]
\includegraphics[height=3cm]{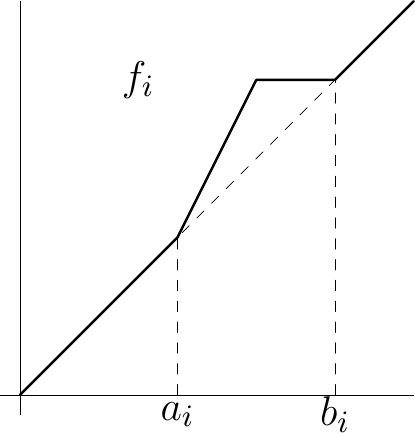}
\caption{The map $f_i$.}
\end{figure}

As neither $0$ nor $1$ belongs to $K$, there is a rational number $\epsilon>0$ such that $K\subseteq [\epsilon,1-\epsilon]$. Let us define $f:[0,1]\to \R$ by $f(x)=x$ on $[\epsilon,1-\epsilon]$, $f(x)=2x-(1-\epsilon)$ on $[1-\epsilon,1]$ and $f(x)=\epsilon$ on $[0,\epsilon]$.

We then define the map $T:[0,1]\to \R$ by $T(x)=\frac{f}{2}+\sum_{i\geq 2} 2^{-i}f_i$. $T$ is computable and non-decreasing, and $T(x)>x$ if and only if $x\in [0,1]\setminus K$. As $T(1)=1+\frac{\epsilon}{2}=1+T(0)$, we can take the quotient of $T$ modulo $1$.

\begin{figure}[h]
\includegraphics[height=3cm]{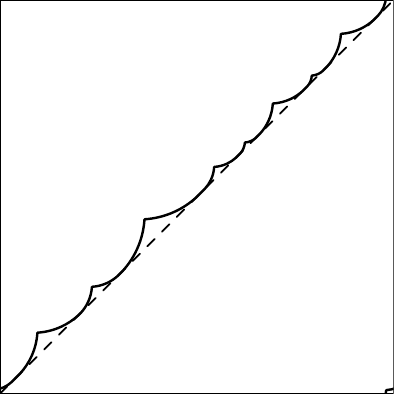}
\caption{The map $T$.}
\end{figure}

\begin{proposition}
$W=U\cup [0,\epsilon) \cup (1-\epsilon , 1] $ is a strictly invariant set: $T^{-1}W=W$.
\end{proposition}

\begin{proof}
If $x\notin W$ then $T(x)=x\notin W$.

If $x\in W$ then $T(x)\in W$. Indeed, if $T(x)\notin W$, $T(x)$ is a fixed point so $T$ is constant on $[x,T(x)]$ ($T$ is non-decreasing). Let $q$ be any rational number in $(x,T(x))$: $T(x)=T(q)$ is then computable, but does not belong to $W$: impossible.
\end{proof}

\begin{proposition}\label{bomb}
The map $T$ is computable but has no computable invariant probability measure.
\end{proposition}

Let $x\in [0,1]$: the trajectory of $x$ is "non-decreasing"  and converges to the first point above $x$ which is not in $U$, $\inf ([x,1]\setminus U)$ or to $\min(K)$ if $x>\sup(K)$. More precisely, there are two cases: (i) if $x\notin U$ then $x$ is a fixed point (unstable on the right), (ii) if $x\in U$ then the trajectory of $x$ converges to a lower semi-computable fixed point (non-computable, as it does not belong to $U$).

\begin{lemma}
Let $\mu$ be an invariant probability measure: then $\mu(K^c)=0$.
\end{lemma}

\begin{proof}
Obviously $\mu (0)=0 $ because $0$ is not periodic. Let $(a,b)=(a_i,b_i)$ be an interval from the description of $U$. Since $T^n(a)$ and $T^n(b)$ tends to some non computable $\alpha$ (and then are not stationary, as they are computable), the interval $(a,b)$ is wandering. Hence, by Poincar\'e recurrence theorem it has null measure.
\end{proof}

\begin{proof}(of proposition \ref{bomb})
We can conclude: let $\mu$ be a computable invariant probability measure: by the above lemma its support is then included in the complement of $W$. But the support of a computable probability measure always contains computable points (see proposition \ref{compoint}) : contradiction.
\end{proof}

Actually, the set of invariant measures is exactly the set of measures which give null weight to $W$. It is easy to see that in the above system the set of invariant measures is a convex recursive compact set. Indeed, the function $\mu\to\mu(W)$ is lower semi-computable, so $\{\mu:\mu(W)>0\}$ is a recursive open set. Its complement is then a recursive compact set, as the whole space of probability measures is a recursive compact set.
The above example hence shows an example of a convex, and recursive compact set whose extremal points are not computable.

We end remarking that with a different construction of the various $f_i$ it is possible to give  also a smooth  system having the same properties as the examples in this section.


\begin{thebibliography}{GHR07}

\bibitem[AGS]{AGS} L. Ambrosio, N. Gigli, G. Savare. \emph{Gradient
flows: in metric spaces and in the space of probability measures}, Birkhauser
Zurich 2005 ebnisse der Mathematik und ihrer Grenzgebiete.


\bibitem[BY06]{BY06} M. Braverman, M. Yampolsky. \emph{Non-computable Julia sets}, Journ. Amer. Math. Soc., 19 (2006), 551-578.

\bibitem[BBY07]{BBY07} I. Binder, M. Braverman, M. Yampolsky. \emph{Filled Julia sets with empty interior are computable}, Journal FoCM, 7(2007), 405-416.

\bibitem[BBY06]{BBY06} I. Binder, M. Braverman, M. Yampolsky. \emph{On computational complexity of Siegel Julia sets}, Commun. Math. Phys., 264, 317-334(2006)

\bibitem[Bla94]{Bla94} M. L. Blank. \emph{Pathologies generated by round-off in dynamical systems}. Physica D. (1994) vol 78, no 1-2, pp: 93--114.

\bibitem[Bla89]{Bla89} M. L. Blank. \emph{Small perturbations of chaotic dynamical systems}. Russian Mathematical Surveys. (1989) vol 44, no 6, pp:1--33.

\bibitem[BP03]{BraPre03}
Vasco Brattka and Gero Presser.
\newblock Computability on subsets of metric spaces.
\newblock {\em Theoretical Computer Science}, 305(1-3):43--76, 2003.

\bibitem[B83]{B83} A. A. Brudno (1983) \emph{Entropy and the complexity of
the trajectories of a dynamical system.} Trans. Mosc. Math. Soc. \textbf{44}
127--151.

\bibitem[CCS]{CCS} J-R. Chazottes, P. Collet and B. Schmitt.
\emph{Statistical consequences of the Devroye inequality for processes. Applications to a class of non-uniformly hyperbolic dynamical systems} 2005 Nonlinearity 18 2341-2364.  


\bibitem[Del99]{Del99} M. Dellnitz, A. Hohmann. \emph{On the Approximation of Complicated Dynamical Behavior}. SIAM Journal on Numerical Analysis. 1999, vol. 36, no2, pp. 491-515.

\bibitem[Del97]{Del97} M. Dellnitz, A. Hohmann. \emph{A subdivision algorithm for the computation of unstable manifolds and global attractors}. Numerische Mathematik. 1997, vol. 75, no3, pp. 293-317.

\bibitem[Din94]{Din94} J. Ding, A. Zhou. \emph{The projection method for computing multidimensional absolutely continuous invariant measures}. Journal of Statistical Physics.(1994) vol 77, 3-4, pp: 899-908.

\bibitem[Din93]{Din93} J. Ding, Q. Du, T. Y. Li. \emph{High order approximation of the Frobenius-Perron operator}. Applied Mathematics and Computation. (1993) vol 53, pp: 151 - 171. 

\bibitem[G93]{G93} P. G\'{a}cs. \emph{Lectures notes on descriptional
complexity and randomness.} Boston University (1993) 1--67.

\bibitem[Ga00]{Ga00} S. Galatolo. \emph{Orbit complexity by computable structures} Nonlinearity \textbf{13}, 1531-1546 (2000).

\bibitem[GHR09a]{GHR08}
P. G{\'a}cs, M. Hoyrup, C. Rojas.
\newblock Randomness on computable probability spaces - a dynamical point of
  view.
\newblock In Susanne Albers and Jean-Yves Marion, editors, {\em 26th
  International Symposium on Theoretical Aspects of Computer Science (STACS
  2009)}, pages 469--480, Dagstuhl, Germany, 2009. Schloss Dagstuhl -
  Leibniz-Zentrum fuer Informatik, Germany.



\bibitem[GHR09b]{GHR07}
S. Galatolo, M. Hoyrup, C. Rojas.
\newblock A constructive {B}orel-{C}antelli lemma. {C}onstructing orbits with
  required statistical properties.
\newblock {\em Theoretical Computer Science}, 2009.
\newblock In Press. doi:10.1016/j.tcs.2009.02.010  (Available on ArXiv).


\bibitem[GHR09c]{GHR09c}
S. Galatolo, M. Hoyrup, C. Rojas.
\newblock Effective symbolic dynamics, random points, statistical behavior, complexity and entropy
\newblock  Available from arXiv. http://arxiv.org/abs/0801.0209


\bibitem[GP09]{GP09} S. Galatolo, M. J. Pacifico. \emph{Lorenz like flows: exponential decay of correlations for the Poincar\'e map, logarithm law, quantitative recurrence. } http://front.math.ucdavis.edu/0901.0574

\bibitem[H05]{H05} P. Hertling \emph{Is the Mandelbrot set computable},  Math. Logic Quart, 2005, vol 51, pp. 5--18.


\bibitem[HR09]{HR07} M. Hoyrup, C. Rojas. \emph{Computability of probability
measures and Martin-L\"of randomness over metric spaces.} Information and Computation, in press. arXiv:0709.0907
(2009)

\bibitem[H95]{H} B. Hunt. \emph{Estimating invariant measures and
Lyapunov exponents} \ preprint - http://www.chaos.umd.edu/\symbol{126} bhunt/research/eimale.pdf (1995)

\bibitem[I03]{I03} S. Isola. \emph{On systems with finite ergodic degree.}
Far East Journal of Dynamical Systems 5(2003), 1-62

\bibitem[KMY98]{KMY} M. Keane, R. Murray, L. S. Young. \emph{Computing
invariant measures for expanding circle maps } Nonlinearity 11 27-46 (1998)

\bibitem[Kif86]{Kif86} Y. Kifer. \emph{General random perturbations of hyperbolic and expanding transformations}. Journal d'Analyse Math\'ematique. (1986) vol 47, 1, pp: 111-150.

\bibitem[LM08]{LM08} L.  Bienvenu, W. Merkle. \emph{Effective Randomness for Computable Probability Measures} Electr. Notes Theor. Comput. Sci.
vol {167}, (2007) p. {117-130}.


\bibitem[L01]{L} C. Liverani. \emph{Rigorous numerical investigations of the
statistical properties of piecewise expanding maps-A feasibility study.}
Nonlinearity 14, n.3 pp. 463-490, (2001).

\bibitem[Pal00]{Pal00} K. Palmer. \emph{Shadowing in Dynamical Systems - Theory and Applications}. Mathematics and Applications vol 501. Springer (2000).

\bibitem[PJ99]{PJ99} M. Pollicott, O. Jenkinson. \emph{Computing Invariant
Densities and Metric Entropy }Comm. Math. Phys.  211, n. 3 \ (2000) 687-703



\bibitem[Rog87]{Rog87} H. Rogers. \emph{Theory of recursive functions and effective computability} MIT Press Cambridge, MA, USA (1987)

\bibitem[Tuc99]{Tuc99} W. Tucker. \emph{The Lorenz attractor exists}, C. R. Acad. Sci. Paris, 1999 vol 328, pp. {1197--1202}.


\bibitem[T36]{T36} A. Turing. \emph{On computable numbers, with an
application to the Entscheidungsproblem }Proc. Lond. Math. Soc. 2, 42
230-265 (1936)


\bibitem[Y02]{Y} L.S. Young. \emph{What are SRB measures, and which dynamical systems have them?} Journal of Statistical Physics, Volume 108, Numbers 5-6, September 2002 , pp. 733-754(22)



\bibitem[V97]{V97} M. Viana. \emph{Stochastic dynamics of deterministic
systems} Brazillian Math. Colloquium 1997, IMPA.


\bibitem[W00]{Wei00} K. Weihrauch. \emph{Computable Analysis}, Springer, 2000.M. 
\end{thebibliography}
\end{document}